\documentclass{conm-p-l}
\usepackage{amsmath,amsthm,latexsym,amssymb}
\usepackage[all]{xy}

\newtheorem{theorem}{Theorem}[section]

\newtheorem{proposition}[theorem]{Proposition}
\newtheorem{corollary}[theorem]{Corollary}

\theoremstyle{definition}
\newtheorem{definition}[theorem]{Definition}
\newtheorem{example}[theorem]{Example}

\theoremstyle{remark}
\newtheorem{remark}[theorem]{Remark}

\numberwithin{equation}{section}

\newdir{ >}{{}*!/-5pt/@{>}}

\newcommand\fib{\ar @{->>} [r]} 
 
\newcommand\we{\ar [r]^{\simeq}}
\newcommand{\hh}{\operatorname H}
\newcommand{\Om}{\Omega}
\newcommand{\vp}{\varphi}
\newcommand{\cdga}{\mathbf{CDGA}_{\mathbb Q}}
\newcommand{\del}{\partial}
\newcommand{\apl}{\mathcal A_  {\text{PL}}}
\newcommand\im{\operatorname {Im}}
\newcommand\fdim{\operatorname {fdim}}
\newcommand\cat{\operatorname{cat}}
\newcommand\mcat{\operatorname{Mcat}}
\newcommand{\ob}{\operatorname{Ob}}

\begin{document}

\title[Rational Homotopy Theory]{Rational Homotopy Theory:\\ A Brief Introduction}

\author{Kathryn Hess}
\address{EPFL SB IGAT, B\^atiment BCH, CH-1015 Lausanne, Switzerland}
\thanks{Work completed during a visit to the Institut Mittag-Leffler (Djursholm,
Sweden)}
\email{kathryn.hess@epfl.ch}

\date{April 10, 2006}

\begin{abstract}
These notes contain a brief introduction to rational homotopy theory: its
model category foundations, the Sullivan model and interactions with the 
theory of local commutative rings.
\end{abstract}

\maketitle
\tableofcontents

\section*{Introduction}
This overview of rational homotopy theory consists of an extended version of lecture notes from a minicourse based primarily on the encyclopedic text \cite {fht} of F\'elix, Halperin and Thomas.  With only three hours to devote to such a broad and rich subject, it was difficult to choose among the numerous possible topics to present.  Based on the subjects covered in the first week of this summer school, I decided that the goal of this course should be to establish carefully the foundations of rational homotopy theory, then to treat more superficially one of its most important tools, the Sullivan model.  Finally, I provided a brief summary of the extremely fruitful interactions between rational homotopy theory and local algebra, in the spirit of the summer school theme ``Interactions between Homotopy Theory and Algebra.''  I hoped to motivate the students to delve more deeply into the subject themselves, while providing them with a solid enough background to do so with relative ease.

As these lecture notes do not constitute a history of rational homotopy theory, I have chosen to refer the reader to \cite {fht}, instead of to the original papers, for the proofs of almost all of the results cited, at least in Sections 1 and 2.  The reader interested in proper attributions will find them in \cite {fht} or \cite {hess}.

The author would like to thank Chris Allday, Luchezar Avramov, Srikanth Iyengar, Mike Mandell and Jonathan Scott, as well as the anonymous referee, for their helpful comments on an earlier version of this article.

\subsection*{Basic notation and terminology}
We assume in this chapter that the reader is familiar with the elements of the theories of simplicial sets and of model categories.  As references we recommend \cite{dwyspal} or \cite{hovey} or the chapter of these lecture notes by Paul Goerss \cite{goerss}.

In this chapter, $\mathbf{sSet}$ and $\mathbf{Top}$ are the categories of simplicial sets and of topological spaces, respectively.  Furthermore, $|\cdot|:\mathbf{sSet}\to \mathbf{Top}$ denotes the geometric realization functor, while $S_\bullet:\mathbf{Top}\to \mathbf{sSet}$ denotes its right adjoint, the singular simplices functor.

If $K$ is a simplicial set, then $C_*(K)$ and $C^*(K)$ denote its normalized chain and cochain complexes, respectively.  If $X$ is a topological space, then  the singular chains and cochains on $X$ are $S_*(X):=C_*\big (S_\bullet (X)\big)$ and $S^*(X):=C^*\big (S_\bullet (X)\big)$.

A morphism of (co)chain complexes inducing an isomorphism in (co)homology is called a \emph{quasi-isomorphim} and denoted $\xrightarrow{\simeq}$.

A graded vector space is said to be of \emph{finite type} if it is finite dimensional in each degree. A CW-complex is of finite type if it has finite number of cells in each dimension.

Given a category $\mathbf{C}$ and two objects $A$ and $B$ in $\mathbf C$, we write $\mathbf C(A,B)$ for the class of morphisms with source $A$ and target $B$.

\section{Foundations}

For the sake of simplicity, we work throughout these notes only with simply connected spaces.  Many of the results presented hold for connected, nilpotent spaces as well.

\subsection{Rationalization and rational homotopy type}

Let $\widetilde\hh _{*}$  denote reduced homology.

\begin{definition} A simply connected space $X$ is \emph{rational} if the following, equivalent conditions are satisfied.
\begin{enumerate}
\item $\pi_*X$ is a $\mathbb Q$-vector space.
\item $\widetilde\hh_*(X;\mathbb Z)$ is a $\mathbb Q$-vector space.
\item $\widetilde\hh_*(\Om X;\mathbb Z)$ is a $\mathbb Q$-vector space.
\end{enumerate}
\end{definition}

To prove the equivalence of these conditions, one begins by observing that $\hh_*(K(\mathbb Q,1);\mathbb F_p)\cong  \hh_*(\text{pt.};\mathbb F_p)$ for all primes $p$. An inductive Serre spectral sequence argument then shows that $\hh_*(K(\mathbb Q,n);\mathbb F_p)\cong  \hh_*(\text{pt.};\mathbb F_p)$ for all primes $p$ and all $n\geq 1$.  The equivalence of conditions (1) and (2) for an arbitrary $X$ then follows from an inductive argument on the Postnikov tower of $X$. On the other hand, the equivalence between conditions (2) and (3) can be easily verified by a Serre spectral sequence argument.

\begin{example} For any $n\geq 2$, let $\iota_{n,k}$ denote the homotopy class of the inclusion of $S^n$ as the $k^{\text{th}}$ summand $S^n_k$  of $\bigvee _{k\geq 1}S^n_k$.  The \emph{rational $n$-sphere} is defined to be the complex
$$(S^n)_0 :=\bigg(\bigvee _{k\geq 1}S^n_k\bigg)\bigcup \bigg(\coprod _{k\geq 2}D^{n+1}_k\bigg),$$
where $D^{n+1}_{k+1}$ is attached to $S^n_k\vee S^n_{k+1}$ via a representative
$S^n\to S^n_k\vee S^n_{k+1}$ of $\iota_{n,k}-(k+1)\iota _{n,k+1}$.  The \emph{rational $n$-disk} is then 
$$(D^{n+1})_0:=(S^n)_0\times I/(S^n)_0\times \{0\}.$$

Let 
$$X(r):=\bigg(\bigvee _{1\leq k\leq r}S^n_k\bigg)\bigcup \bigg(\coprod _{2\leq k\leq r-1}D^{n+1}_k\bigg).$$
It is clear that for all $r$, $S^n_r$ is a strong deformation retract of $X(r)$, which implies that $\hh_kX(R)=0$ if $k\not= 0,n$. Furthermore, the homomorphism induced in reduced homology by the inclusion $X(r)\hookrightarrow X(r+1)$ is multiplication by $r+1$.  Since homology commutes with direct limits and $(S^n)_0=\underset {\to}{\lim}\, X(r)$,
$$\widetilde \hh_*\big ((S^n)_0;\mathbb Z)=\left\{\begin{array} {r@{\quad:\quad}l}
					\mathbb Q& k=n      \\ 
					0&\text{else.}
\end{array}\right.$$ 
\end{example}

\begin{definition} A pair of spaces $(X,A)$ is a \emph{relative $\text{CW}_0$-complex} if $X=\bigcup_{n\geq 1}X(n)$, where
\begin{enumerate}
\item  $X(1)=A$;
\item for all $n\geq 1$, there is a pushout
$$\xymatrix{
\coprod _{\alpha \in J_n} (S^n)_0\ar [d]_{\text{incl.}}\ar[r]^{\coprod _\alpha f_\alpha}&X(n)\ar [d]\\
\coprod _{\alpha \in J_n} (D^{n+1})_0 \ar[r]&X(n+1);}$$
\item $X$ has the weak topology, i.e., $U\subset X$ is open in $X$ if and only if $U\cap X(n)$ is open in $X(n)$ for all $n$.
\end{enumerate}
\end{definition}

The pairs $\big((S^n)_0, S^n\big )$ and $\big((D^{n+1})_0, D^{n+1}\big )$ are the fundamental examples of relative $\text{CW}_0$-complexes.

Note that if $A$ is a rational space and $(X,A)$ is a relative $\text{CW}_0$-complex, then $X$ is rational as well.

\begin{definition} Let $X$ be a simply connected space. A continuous map $\ell: X\to Y$ is a \emph{rationalization} of $X$ if $Y$ is simply connected and rational and
$$\pi_*\ell\otimes \mathbb Q:\pi_*X\otimes \mathbb Q\longrightarrow  \pi _*Y\otimes \mathbb Q\cong \pi _*Y$$
is an isomorphism.
\end{definition}

Observe that a map $\ell: X\to Y$ of simply connected spaces is a rationalization if and only if $\hh _*(\ell;\mathbb Q)$ is an isomorphism and $Y$ is rational.

The inclusions of $S^n$ into $(S^n)_0$ and of $D^{n+1}$ into $(D^{n+1})_0$ are rationalizations.  The rationalization of an arbitrary simply connected space, as constructed in the next theorem, generalizes these fundamental examples.

\begin{theorem}\label{thm:ratlztn}Let $X$ be a simply connected space.  There exists a relative CW-complex $(X_0,X)$ with no zero-cells and no one-cells such that the inclusion $j:X\to X_0$ is a rationalization.  Furthermore, if $Y$ is a simply connected, rational space, then any continuous map $f:X\to Y$ can be extended over $X_0$, i.e., there is a continuous map $g:X_0\to Y$, which is unique up to homotopy, such that 
$$\xymatrix{
X\ar[rr]^f\ar[dr]_j&&Y\\ &X_0\ar[ur]_g}$$
commutes.
\end{theorem}

\begin{proof} We provide only a brief sketch of the proof.  We can restrict to the case where $X$ is a $1$-reduced CW-complex. The rationalization $X_0$ can then be constructed as a $\text{CW}_0$-complex with rational $n$-cells in bijection with the $n$-cells of $X$, for all $n$.  The attaching maps of $X_0$ are obtained by rationalizing the attaching maps of $X$.  The complete proof can be found in \cite[Theorem 9.7]{fht}.
\end{proof}

Continuing in the same vein, one can show that such a cellular rationalization is unique up to homotopy equivalence, relative to $X$.

Given a continous map $\vp: X\to Y$ between simply connected spaces, we let $\vp_0:X_0\to Y_0$ denote the induced map between their rationalizations, the existence and uniqueness (up to homotopy) of which are guaranteed by Theorem \ref{thm:ratlztn}.

\begin{definition} The \emph{rational homotopy type} of a simply connected space $X$ is the weak homotopy type of its rationalization $X_0$.
\end{definition}

\begin{definition} A continuous map $\vp:X\to Y$ between simply connected spaces is a \emph {rational homotopy equivalence} if the following, equivalent conditions are satisfied.
\begin{enumerate}
\item $\pi_*(\vp)\otimes \mathbb Q$ is an isomorphism.
\item $\hh_*(\vp;\mathbb Q)$ is an isomorphism.
\item $\hh^*(\vp;\mathbb Q)$ is an isomorphism.
\item $\vp_0:X_0\to Y_0$ is a weak homotopy equivalence.
\end{enumerate}
\end{definition}

To simplify computations, it is common in rational homotopy theory to restrict to the class of spaces defined by the following proposition, the proof of which is in \cite[Theorem 9.11]{fht}.

\begin{proposition}\label{prop:fin-typ} For any simply connected space $X$, there is a CW-complex $Z$ and a rational homotopy equivalence $\vp:Z\to X$ such that
\begin{enumerate}
\item $\hh_*(X;\mathbb Q)$ is of finite type if and only if $Z$ is of finite type; and
\item if $\dim _{\mathbb Q} \hh_*(X; \mathbb Q)<\infty$, then $\hh_*(X;\mathbb Q)=\hh_{\leq N}(X;\mathbb Q)$ if and only if $Z$ is a finite CW-complex of dimension at most $N$.
\end{enumerate}
\end{proposition}

\begin{definition} A simply connected space $X$ is of \emph{finite rational type} if condition (1) of Proposition \ref{prop:fin-typ} is satisfied.
\end{definition}

We can now finally specify clearly the subject presented in these notes.

\begin{quote}
\textbf{Rational homotopy theory is the study of rational homotopy types of spaces and of the properties of spaces and maps that are invariant under rational homotopy equivalence.}
\end{quote}

For further information on rationalization, the reader is refered to Section 9 of  \cite{fht}.

\subsection{The passage to commutative cochain algebras}

We show in this section that the category of rational homotopy types of simply connected, finite-type spaces and of homotopy classes of maps between their representatives is equivalent to an appropriately defined homotopy category of commutative differential graded algebras over $\mathbb Q$.

\subsubsection*{The algebraic category and its homotopy structure}  We begin by a rather careful introduction to the algebraic category in which the Sullivan model of a topological space lives.

A \emph{commutative differential graded algebra (CDGA)} over $\mathbb Q$ is a commutative monoid in the category of non-negatively graded, rational cochain complexes (cf., \cite {goerss}).  In other words, a CDGA is a cochain complex $(A^*,d)$ over $\mathbb Q$, endowed with cochain maps
$$\eta: \mathbb Q\longrightarrow  (A^*,d),$$
called the \emph{unit}, and
$$\mu:(A^*,d)\otimes_{\mathbb Q} (A^*,d)\longrightarrow  (A^*,d):a\otimes b\mapsto a\cdot b,$$
called the \emph{product}, such that
\begin{enumerate}
\item $\mu$ is graded commutative, i.e., if $a\in A^p$ and $b\in A^q$, then $a\cdot b=(-1)^{pq} b\cdot a$;
\item $\mu$ is associative; and
\item $\mu(\eta\otimes Id_A)=Id_A=\mu(Id_A\otimes \eta)$.
\end{enumerate}
Observe that since $\mu $ is a morphism of chain complexes, the differential on a CDGA satisfies the Leibniz rule, i.e., if $a\in A^p$ and $b\in A^q$, then
$$d(a\cdot b)=da\cdot b +(-1)^p a\cdot db.$$

Let $r\in \mathbb N$.  A CDGA $A$ is \emph{$r$-connected} if $A^0=\mathbb Q$ and $A^k=0$ for all $0<k<r+1$.

A morphism of CDGA's $f:(A^*,d,\mu, \eta)\to(\bar A^*,\bar d,\bar \mu, \bar \eta)$ is a cochain map such that $f\mu=\bar\mu(f\otimes f)$ and $f\eta=\bar\eta$.  The category of CDGA's over $\mathbb Q$ and their morphims  is denoted $\cdga$.

To simplify notation, we frequently write either $A$ or $(A,d)$ to denote $(A^*,d,\mu, \eta)$. Furthermore, henceforth in these notes, the notation $\otimes$ means $\otimes_{\mathbb Q}$.

In rational homotopy theory, CDGA's with free underlying commutative, graded algebra play an essential role. Given a non-negatively graded vector space $V=\oplus_{i\geq 0}V^i$, let $\Lambda V$ denote the free, commutative, graded algebra generated by $V$, i.e.,
$$\Lambda V=S[V^{\text{even}}] \otimes E[V^{\text{odd}}],$$
the tensor product of the symmetric algebra on the vectors of even degree and of the exterior algebra on the vectors of odd degree. Given a basis $\{v_j\mid j\in J\}$ of $V$, we often write $\Lambda (v_j)_{j\in J}$ for $\Lambda V$.  We also write $\Lambda ^nV$ to denote the set of elements of $\Lambda V$ of wordlength $n$.

A homomorphism of commutative, graded algebras $\vp :\Lambda V\to A$ is determined by its restriction to $V$, as is any derivation of commutative, graded algebras $\delta:\Lambda V\to A$. In particular, the differential $d$ of a CDGA $(\Lambda V,d)$ is determined by its restriction to $V$.

More generally, the following class of CDGA morphisms is particularly important in rational homotopy theory.

\begin{definition} A \emph{relative Sullivan algebra} consists of an inclusion of CDGA's $(A,d)\to (A\otimes\Lambda V,d)$ such that $V$ has a basis $\{ v_\alpha\mid \alpha \in J\}$, where $J$ is a well-ordered set, such that $dv_\beta\in A\otimes\Lambda V_{<\beta}$ for all $\beta\in J$, where $V_\beta$ is the span of $\{v_\alpha \mid \alpha<\beta\}$. A relative Sullivan algebra is \emph{minimal} if
$$\alpha<\beta\Longrightarrow  \deg v_\alpha\leq\deg v_\beta.$$
A (minimal) relative Sullivan algebra $(\Lambda V, d)$ extending $(A,d)=(\mathbb Q, 0)$ is called a \emph{(minimal) Sullivan algebra}.
\end{definition}

\begin{remark} If $V^0=0=V^1$, then $(\Lambda V,d)$ is always a Sullivan algebra and is minimal if and only if $dV\subset \Lambda ^{\geq 2} V$, the subspace of words of length at least two.
\end{remark}

\begin{example} The CDGA $(\Lambda (x,y,z),d)$, where $x$, $y$ and $z$ are all of degree $1$ and 
$dx=yz$, $dy=zx$ and $dz=yx$, is an example of a CDGA with free underlying graded algebra that is not a Sullivan algebra.
\end{example}

Recall from Example 1.7 and (a slight modification of) Example 3.4(1) in \cite{goerss} that the category $\mathbf{Ch}^*(\mathbb Q)$ of non-negatively graded cochain complexes over $\mathbb Q$ admits a cofibrantly generated model category structure in which
\begin{enumerate}
\item weak equivalences are quasi-isomorphisms;
\item fibrations are degreewise surjections; and
\item cofibrations are degreewise injections, in positive degrees.
\end{enumerate}

The set of generating acyclic cofibrations is
$$\mathcal J=\{j_n:0\to D(n)\mid n\geq 1\}$$
where
$$D(n)^k=\left\{\begin{array} {r@{\quad:\quad}l}
		\mathbb Q& k=n-1,n\\
		0&\text{else}
	\end{array}\right.
$$
and $d:D(n)^{n-1}\to D(n)^n$ is the identity map.  The set $\mathcal I$ of generating cofibrations consists of the inclusions
$$i_n:S(n)\to D(n)\quad\text{for}\quad n\geq 1\qquad\text{and}\qquad i_0:0\to S(0), i_{0}':S(0)\to 0,$$
where 
$$S(n)^k=\left\{\begin{array} {r@{\quad:\quad}l}
		\mathbb Q& k=n\\
		0&\text{else.}
	\end{array}\right.
$$

Consider the pair of adjoint functors
$$\Lambda:\mathbf{Ch}^*(\mathbb Q) \rightleftarrows \cdga:U,$$
where $\Lambda $ is the ``free commutative cochain algebra'' functor satisfying $\Lambda (C,d)=(\Lambda C,\bar d)$, where $\bar d$ is the derivation extending $d$, and $U$ is the forgetful functor.  It is not difficult, as indicated in Example 3.7 in \cite {goerss}, to show that this adjoint pair satisfies the hypotheses of Theorem 3.6 of \cite{goerss}.

There is thus a cofibrantly generated model structure on $\cdga$, with generating cofibrations and acyclic cofibrations
$$\mathcal I=\{ \Lambda i_n, \Lambda i_{0}'\mid n\geq 0\}\quad\text{and}\quad \mathcal J=\{\Lambda j_n\mid n\geq 1\},$$
where $\Lambda (0):=\mathbb Q$.  Let $\mathcal I-cell$ denote the smallest class of morphisms in $\cdga$ that contains $\mathcal I$ and that is closed under coproducts, cobase change and sequential colimits.  It is easy to see that $\mathcal I-cell$ is exactly the class of relative Sullivan algebras.  In this model structure on $\cdga$, weak equivalences are quasi-isomorphisms, and fibrations are degreewise surjections. Cofibrations are retracts of relative Sullivan algebras.  All CDGA's are fibrant, and the Sullivan algebras are the cofibrant CDGA's.

The next result, which is the Lifting Axiom of model categories in the specific case of $\cdga$, is an important tool in rational homotopy theory.

\begin{proposition}[The Lifting Lemma]\label{prop:lift} Let
$$\xymatrix{
(A,d)\ar[d]_i\ar[r]^f&(B,d)\ar[d]_p\\
(A\otimes \Lambda V,D)\ar[r]^g&(C,d)}
$$
be a commuting diagram in $\cdga$, where $i$ is a relative Sullivan algebra and $p$ is a surjection.  If $i$ or $p$ is a quasi-isomorphism, then $g$ lifts through $p$ to an extension of $f$, i.e., there exists a CDGA map $h:(A\otimes \Lambda V,D)\to (B,d)$ such that $hi=f$ and $ph=g$.  Furthermore, any two lifts are homotopic relative to $(A,d)$.
\end{proposition}

We can describe homotopy of CDGA morphisms with source a Sullivan algebra in terms of the following path objets.  Let $I$ denote the CDGA $D(1)$, where the generators of degrees $0$ and $1$ are called $t$ and $y$, respectively.  Let $\varepsilon _0:I\to\mathbb Q$ and $\varepsilon _1:I\to\mathbb Q$ denote the augmentations specified by $\varepsilon_0(t)=0$ and $\varepsilon _1(t)=1$.

\begin{proposition} Let $(\Lambda V,d)$ be a Sullivan algebra.  Two CDGA morphisms $f,g:(\Lambda V,d)\to (A,d)$ are homotopic if and only if there is a CDGA morphism $H:(\Lambda V,d)\to(A,d)\otimes I$ such that $(Id_A\otimes \varepsilon _0)H=f$ and $(Id_A\otimes \varepsilon _1)H=g$.
\end{proposition}

A careful, degree-by-degree version of the proof of the Small Object Argument (Theorem 3.5 in \cite{goerss}) establishes the following useful result.

\begin{proposition}  Any morphism $f:(A,d)\to(B,d)$ in $\cdga$  can be factored as
$$\xymatrix{
&(A\otimes \Lambda U,D)\ar[dr]^p\\
(A,d)\ar[ur]^i_{\simeq}\ar[rr]^f\ar[dr]_j&&(B,d)\\
&(A\otimes\Lambda V,D)\ar[ur]_q^\simeq}
$$
where $i$ and $j$ are relative Sullivan algebras, and $p$ and $q$ are surjections.
\end{proposition}

If we are willing to sacrifice surjectivity of $q$, we can obtain minimality of $j$, again via a degree-by-degree construction.

\begin{proposition} Any morphism $f:(A,d)\to(B,d)$ in $\cdga$  can be factored as
$$\xymatrix{
(A,d)\ar[rr]^f\ar[dr]_\iota&&(B,d)\\
&(A\otimes\Lambda W,D)\ar[ur]_\vp^\simeq}
$$
where $\iota $ is a relative Sullivan algebra and $\vp$ is a quasi-isomorphism.  In particular, if $\hh^0A\cong \mathbb Q$, $\hh^1A=0$ and $\hh^*B$ is of finite type, then $W$ is of finite type and $W=W^{\geq 2}$.
\end{proposition}

\begin{definition} The quasi-isomorphism $\vp: (A\otimes \Lambda W,D)\xrightarrow {\simeq}(B,d)$ is a \emph{relative Sullivan minimal model} of $f:A\to B$.  A \emph{Sullivan minimal model} of the CDGA $(B,d)$ is a relative Sullivan minimal model $\vp: (\Lambda W,d)\xrightarrow  {\simeq}(B,d)$ of the unit map $\eta:(\mathbb Q,0)\to (B,d)$.
\end{definition}

It is very convenient to know that (relative) Sullivan minimal models are unique up to isomorphism, which is an immediate consequence of the following proposition.

\begin{proposition}  Suppose that
$$\xymatrix{
(A,d)\ar[d]_i\ar[r]^f_\cong&(B,d)\ar[d]_j\\
(A\otimes \Lambda V,D)\ar[r]^{\hat f}_\simeq &(B\otimes \Lambda W,D)}
$$
is a commuting diagram in $\cdga$, where $i$ and $j$ are minimal relative Sullivan algebras, $f$ is an isomorphism and $\hat f$ is a quasi-isomorphism.  Then $\hat f$ is also an isomorphism.
\end{proposition}

The proof of this proposition reduces to showing that if a CDGA endomorphism of a minimal relative Sullivan algebra $(B\otimes \Lambda W,D)$ that fixes $B$ is homotopic to the identity, then it is equal to the identity.

\subsubsection*{The functors} We now explain the passage from topology to algebra, starting with the relationship between simplicial sets and CDGA's.

\begin{definition} The \emph{algebra of polynomial differential forms}, denoted $\mathfrak A_\bullet^*$, is the simplicial CDGA given by
$$\mathfrak A_n^*=\big(\Lambda (t_0,...,t_n;y_0,...,y_n)/J_n,d\big),$$
where $\deg t_i=0$ and $dt_i=y_i$ for all $i$ and $J_n$ is the ideal generated by $\{1-\sum _{i=0}^n t_i, \sum _{j=0}^n y_j\}$.  The faces and degeneracies are specified by
$$\del _i:\mathfrak A_n^* \longrightarrow  \mathfrak A_{n-1}^*:t_k \mapsto \left\{\begin{array} {r@{\quad:\quad}l}
		t_k&k<i\\
		0&k=i\\
		t_{k-1}&k>i
\end{array}\right. $$
and
$$s _i:\mathfrak A_n^* \longrightarrow  \mathfrak A_{n+1}^*:t_k \mapsto \left\{\begin{array} {r@{\quad:\quad}l}
		t_k&k<i\\
		t_k+t_{k+1}&k=i\\
		t_{k+1}&k>i.
\end{array}\right. $$
\end{definition}

The terminology used in this definition is justified by the following observation.  Let $\Om _{\text{DR}}(\Delta ^n)$ be the cochain algebra of smooth forms on $\Delta^n$, the standard topological $n$-simplex.  Then
$$\Om _{\text{DR}}(\Delta ^n)=C^\infty(\Delta ^n)\otimes_{\mathfrak A^0_n}\mathfrak A^*_n,$$
where the cochain algebra morphisms induced by the topological  face inclusions and degeneracy maps agree with $\del _i$ and $s_i$.

\begin{definition}  Let $\mathcal A^*:\mathbf{sSet}\to\cdga$ be the functor specified by $\mathcal A^*(K)=\mathbf{sSet}(K,\mathfrak A^*_\bullet)$, with product and differential defined objectwise.  For any topological space $X$, let $\mathcal A_  {\text{PL}}(X):=\mathcal A^*\big(S_\bullet(X)\big)$, which we call the CDGA of \emph{piecewise-linear de Rham forms} on $X$.
\end{definition}

Since $\mathcal A^*(K)$ is a commutative algebra for every simplicial set $K$, while $C^*(K;\mathbb Q)$ usually is not, we cannot expect to be able to define a natural quasi-isomorphism of cochain algebras directly from the former to the latter.  However, as explained below, there is a cochain map between them that is close to being an algebra ma

Given $f\in \mathcal A^n(K)$, i.e., $f:K\to \mathfrak A^n_\bullet$, and $x\in K_n$, write
$$f(x)=\hat f(x)dt_1\cdots dt_n,$$
so that $\hat f(x)\in \mathbb Q[t_1,...,t_n]$, the ring of polynomials in $t_1,...,t_n$ with coefficients in $\mathbb Q$.  Define a graded linear map $\oint:\mathcal A^*(K)\to C^*(K;\mathbb Q)$ by
$$\big(\oint f\big)(x)=\int_{\Delta^n}  \hat f(x)dt_1\cdots dt_n.$$

\begin{theorem}[The Polynomial Stokes-De Rham Theorem]\label{thm:poly-stokes} The map $\oint$ is a map of cochain complexes, inducing an isomorphism of algebras in cohomology.
\end{theorem}

The proof of this theorem, which can be found in \cite{bousgug}, Theorem 2.2 and Corollary 3.4, proceeds by methods of acyclic models.

Theorem \ref{thm:poly-stokes} can in fact be strengthened: as proved in Proposition 3.3 in \cite{bousgug}, the cochain map $\oint$ is actually a \emph{strongly homotopy multiplicative map}, in the sense of, e.g., Gugenheim and Munkholm \cite {gugmunk}.

To compare the homotopy theory of CDGA's and of simplicial sets via the functor $\mathcal A^*$, we need $\mathcal A^*$ to be a member of an adjoint pair.  We construct its adjoint as follows.

\begin{definition} Let $\mathcal K_\bullet:\cdga\to \mathbf{sSet}$ be the functor specified by $\mathcal K_\bullet(A)=\cdga(A, \mathfrak A^*_\bullet)$, with faces and degeneracies defined objectwise.
\end{definition}

Let $\mathbf{sCDGA}_{\mathbb Q}$ denote the category of simplicial CDGA's.  Given any CGDA $A$ and any simplicial set $K$, we can form a simplicial CDGA, denoted $A\times K$, where the underlying graded CDGA is the tensor product of copies of $A$, indexed over the simplices of $K$.  Furthermore, there are natural isomorphisms
$$\cdga\big(A,\mathbf{sSet}(K,\mathfrak A^*_\bullet)\big)\cong \mathbf{sCDGA}_{\mathbb Q}(A\times K, \mathfrak A^*_\bullet)\cong \mathbf{sSet}\big(K, \cdga (A, \mathfrak A^*_\bullet)\big),$$
i.e.,
$$\cdga\big (A, \mathcal A^*(K)\big)\cong \mathbf{sSet}\big(K, \mathcal K_\bullet(A)\big).$$
We therefore have an adjoint pair
$$\mathcal A^*:\mathbf{sSet}\rightleftarrows\cdga^{\text{op}}: \mathcal K_\bullet,$$
which Bousfield and Gugenheim proved to be a Quillen pair in Section 8 of \cite {bousgug}.

\begin{definition}  The composite functor
$$\xymatrix{ \cdga\ar[dr]_{\mathcal K_\bullet}\ar[rr]^{<->}&& \mathbf{Top}\\
&\mathbf{sSet}\ar[ur]_{|-|}}$$
is called \emph{spatial realization}.
\end{definition}

Let $(\Lambda V,d)$ be any Sullivan algebra.  Let $\eta:Id\to \mathcal A^*\circ \mathcal K_\bullet$ be the unit of the adjoint pair above, and let $\varepsilon: S_\bullet\circ |-|\to Id$ be the counit of the adjoint pair $(S_\bullet, |-|)$. Consider the commuting diagram
$$\xymatrix{
\mathbb Q\ar[d]\ar[r]& \apl\big(<(\Lambda V,d)>\big)\ar [d] _\simeq^{\apl (\varepsilon _{\mathcal K_\bullet(\Lambda V,d)})}\\
(\Lambda V,d)\ar[r]^(0.4){\eta_{(\Lambda V,d)}}&\mathcal A^*\big(\mathcal K_\bullet(\Lambda V,d)\big).}$$
Since $(\Lambda V,d)$ is a Sullivan algebra and $\apl (\varepsilon _{\mathcal K_\bullet(\Lambda V,d)})$ is a surjective quasi-isomorphism, the Lifting Lemma (Proposition \ref{prop:lift}) can be applied to this diagram, establishing the existence of a CDGA morphism $m_{(\Lambda V,d)}:(\Lambda V,d)\to\apl\big(<(\Lambda V,d)>\big)$, unique up to homotopy, lifting $\eta_{(\Lambda V,d)}$. 

\begin{theorem}\label{thm:real1} If $(\Lambda V,d)$ is a simply connected Sullivan algebra of finite type, then 
\begin{enumerate}
\item $m_{(\Lambda V,d)}:(\Lambda V,d)\to\apl\big(<(\Lambda V,d)>\big)$ is a quasi-isomorphism; and
\item $<(\Lambda V,d)>$ is a simply connected, rational space of finite type, such that there is an isomorphism of graded, rational vector spaces
$$\pi_*\big(<(\Lambda V,d)>\big)\cong \hom _{\mathbb Q}(V,\mathbb Q).$$
\end{enumerate}
\end{theorem}

We refer the reader to Section 17 of \cite{fht} for the details of the proof of this extremely important theorem.

To complete the picture, we need to specify the relationship between spatial realization and homotopy of morphisms.

\begin{theorem}\label{thm:real2} Let $(\Lambda V,d)$ and$(\Lambda W,d)$ be simply connected Sullivan algebras of finite type.
\begin{enumerate}
\item Let $f:(\Lambda V,d)\to (\Lambda W,d)$ be a CDGA morphism.  Then
$$\xymatrix{(\Lambda V,d)\ar[d]_{m_{(\Lambda V,d)}}^\simeq\ar[rr]^f&&(\Lambda W,d)\ar[d]_{m_{(\Lambda WV,d)}}^\simeq\\
\apl\big(<(\Lambda V,d)>\big)\ar[rr]^{\apl\big(<f>\big)}&&\apl\big(<(\Lambda W,d)>\big)}$$
commutes up to homotopy.
\item Two CDGA morphisms $f,g :(\Lambda V,d)\to (\Lambda W,d)$ are homotopic if and only if $<f>$ and $<g>$ are homotopic.
\item Let $\alpha: X\to Y$ be a continuous map between simply connected CW-complexes of finite type.  If there is a homotopy-commutative diagram of CDGA's
$$\xymatrix{
(\Lambda V,d)\ar[d]_{\vp }^\simeq\ar[r]^f&(\Lambda W,d)\ar[d]_{\psi}\\ 
\apl(Y)\ar[r]^{\apl(\alpha)}&\apl(X),}$$
then there is a homotopy-commutative diagram of topological spaces
$$\xymatrix{
X\ar[d]_\beta\ar[r]^\alpha&Y\ar[d]_\gamma\\
<(\Lambda W,d)>\ar[r]^{<f>}&<(\Lambda V,d)>}$$
in which $\pi_*(\beta)\otimes\mathbb Q$ and $\pi_*(\gamma)\otimes\mathbb Q$ are isomorphisms.
\end{enumerate}
\end{theorem}

Again, we refer the reader to Section 17 of \cite {fht} for the proof of this theorem.

\begin{corollary} \begin{enumerate}
\item Rational homotopy types of simply connected spaces of finite rational type are in bijective correspondence with isomorphism classes of minimal Sullivan algebras.  
\item Homotopy classes of continuous maps between simply connected, finite-type rational spaces are in bijective correspondence with homotopy classes of morphisms between simply connected, finite-type Sullivan algebras.
\end{enumerate}
\end{corollary}

We are now ready to introduce one of the very most important tools in rational homotopy theory.

\begin{definition} The \emph{Sullivan minimal model} of a simply connected topological space of finite rational type is the unique (up to isomorphism) Sullivan minimal model of its algebra of piecewise-linear de Rham forms
$$\vp:(\Lambda V,d)\xrightarrow {\simeq}\apl(X).$$
\end{definition}

As a consequence of Theorems \ref{thm:real1} and \ref{thm:real2}, if $\vp:(\Lambda V,d)\xrightarrow {\simeq}\apl(X)$ is a Sullivan minimal model, then there is an isomorphism of graded, rational vector spaces
$$\hom_{\mathbb Q}(V,\mathbb Q)\cong \pi_*(X)\otimes \mathbb Q.$$
In other words, given a Sullivan minimal model of a space, we can read off the nontorsion part of its homotopy groups from the generators of the model.

\section {Sullivan models}

Since the CDGA
$\apl (X)$ is huge and has a complicated product, rational homotopy
theorists prefer to carry out computations with the Sullivan minimal model, which is free as an
algebra, with
only finitely many generators in each dimension if $X$ is of finite rational type.  In this section, we provide a brief overview of the power of the Sullivan model.  We begin by providing a few explicit examples of Sullivan minimal models.  We then explore the relationship between topological fibrations and the Sullivan model.  In particular, we explain the slogan ``the Sullivan model of fiber is the cofiber of the Sullivan model'' and illustrate its application.  A classical and essential numerical homotopy invariant, Lusternik-Schnirelmann category, is our next subject: its elementary properties, how to calculate it using the Sullivan model and its additivity.  Finally, we present the beautiful and striking rational dichotomy of finite CW-complexes, the proof of which depends crucially on Lusternik-Schnirelmann category. 

\subsection {Examples and elementary construction}

As a warmup and an aid to developing the reader's intuition, we calculate a few explicit examples of Sullivan models.  Here, a subscript on a generator always indicates its degree. 

\subsubsection*{Spheres}
The Sullivan model of an odd sphere $S^{2n+1}$ is 
$$\varphi:(\Lambda (x_{2n+1}),0)\longrightarrow \apl (S^{2n+1}),$$ 
where $\varphi (x)$ is any representative of the unique cohomology generator of degree $2n+1$.  Since $\varphi$ is obviously a quasi-isomorphism of CDGA's,  the nontorsion and positive-degree part of $\pi_{*}S^{2n+1}$ is concentrated in degree $2n+1$, where it is of rank $1$.

On the other hand, the Sullivan model of an even sphere $S^{2n}$ is 
$$\varphi: (\Lambda (y_{2n}, z_{4n-1}),d)\longrightarrow \apl (S^{2n}),$$ 
where $dz=y^2$ and $\varphi (y)$ represents the unique cohomology generator of degree $2n$.  Since the square of $\varphi (y)$ must be a boundary, there is an acceptable choice of $\varphi (z)$.  Again, $\varphi$ is clearly a quasi-isomorphism of CDGA's, implying that the nontorsion part of $\pi _{*}S^{2n}$ is concentrated in degrees $2n$ and $4n-1$ and that it is of rank $1$ in each of those degrees.

\subsubsection*{Complex projective spaces}
From the long exact sequences in homotopy of the fibrations
$$S^1\longrightarrow  S^{2n+1}\longrightarrow \mathbb CP^n$$
and 
$$S^1\longrightarrow  S^{\infty}\longrightarrow \mathbb CP^{\infty},$$
and the computation above of $\pi _{*}S^{2n+1}\otimes \mathbb Q$, we conclude that
$$\pi _{*}\mathbb CP^n\otimes \mathbb Q=\mathbb Q\cdot u_{2}\oplus \mathbb Q\cdot x_{2n+1}\text{ and }\pi _{*}\mathbb CP^{\infty}\otimes \mathbb Q=\mathbb Q\cdot u_{2}.$$
Consequently, the Sullivan model of $\mathbb CP^n$ is of the form
$$\varphi: \bigl(\Lambda (u_{2},x_{2n+1}), d\bigr)\longrightarrow \apl (\mathbb CP^n),$$
where $dx=u^{n+1}$, $\varphi (u)$ represents the algebra generator of $H^*(\mathbb CP^n;\mathbb Q)$, which is a truncated polynomial algebra on a generator of degree 2, and $\varphi (x)$ kills its $(n+1)^{st}$ power.  The value of $dx$ is nonzero since $H^*(\mathbb CP^n;\mathbb Q)$ is zero in odd degrees.

The Sullivan model for $\mathbb CP^{\infty}$ is even easier to specify since there can be no nontrivial differential.  It is
$$\varphi:(\Lambda(u),0)\longrightarrow \apl (\mathbb CP^{\infty}),$$
where $\varphi (u)$ represents the algebra generator of $H^*(\mathbb CP^{\infty};\mathbb Q)$, which is a polynomial algebra on a generator of degree 2

\subsubsection*{Products}
Let $(\mathbf B, \otimes)$ and $(\mathbf C, \otimes)$ be monoidal categories.  Recall  that a functor $\mathcal F:(\mathbf B, \otimes)\to (\mathbf C,\otimes)$ is \emph{lax monoidal} if  for all $B,B'\in \ob \mathbf B$, there is a natural morphism $F(B)\otimes F(B')\to F(B\otimes B')$ that is appropriately compatible with the associativity and unit isomorphisms in $(\mathbf B, \otimes)$ and $(\mathbf C, \otimes)$.

It is easy to see that $\apl$ is a lax monoidal functor, via the natural quasi-isomorphism $\alpha _{X,Y}$, defined to be the composite
$$\apl(X)\otimes\apl (Y)\xrightarrow {\apl (p_{1})\otimes \apl (p _{2})}\apl (X\times Y)\otimes \apl (X\times Y)\xrightarrow {\mu} \apl (X\times Y),$$
where $p_{i}$ is projection onto the $i^{th}$ component, and $\mu $ is the product on $\apl (X\times Y)$.

Given Sullivan models $\varphi: (\Lambda V, d)\to\apl (X)$ and $\varphi': (\Lambda V', d')\to\apl(X')$, the Sullivan model of the product space $X\times X'$ is given by
$$(\Lambda V,d)\otimes (\Lambda V', d')\xrightarrow {\varphi\otimes \varphi'}\apl (X)\otimes \apl (X')\xrightarrow {\alpha _{X,X'}}\apl(X\times X').$$

\subsubsection* {Formal spaces}  A space is formal if its rational homotopy is a formal consequence of its rational cohomology, in the sense of the following definition.

\begin{definition}  A rational CDGA $A$ is \emph{formal} if there is a zigzag of quasi-isomorphisms of CDGA's 
$$\xymatrix{A&\bullet\ar[l]_(0.45){\simeq}\ar [r]^(0.35)\simeq& \hh^*(A).}$$  
A space $X$ is \emph{formal} Êif $\apl(X)$ is a formal CDGA.
\end{definition}

From the previous examples, we see that spheres and complex projective spaces are formal.  Furthermore, products of formal spaces are clearly formal.  We can also show that wedges of formal spaces are formal, as follows.

Given a set of augmented CDGA's $\{A_{j}\mid j\in J\}$, let $\prod_{j\in J}^\downarrow A_{j}$ denote their fibered product, i.e., their categorical product in the overcategory $\cdga\negthinspace\downarrow\negthinspace \mathbb Q$. Recall furthermore that for any family of well-pointed spaces $\{X_{j}\mid j\in J\}$ there is a weak equivalence 
$$\bigvee _{j\in J}S_{\bullet } (X_{j})\xrightarrow{\simeq} S_{\bullet } (\bigvee _{j\in J}X_{j}),$$
which induces a quasi-isomorphism
$$\mathcal A^*\big(S_{\bullet } (\bigvee _{j\in J}X_{j})\big)\xrightarrow{\simeq}\mathcal A^*\big(\bigvee _{j\in J}S_{\bullet } (X_{j})\big),$$
since $\mathcal A^*$ is the the left member of a Quillen pair and therefore preserves weak equivalences between cofibrant objects.
It follows that   
\begin{multline*}
\apl (\bigvee _{j\in J}X_{j})=\mathcal A^*\bigl (S_{\bullet}(\bigvee _{j\in J}X_{j})\bigr)\simeq \mathcal A^*\bigl(\bigvee _{j\in J}S_{\bullet}(X_{j})\bigr)\\
\cong\prod _{j\in J}\negthinspace{}^\downarrow\mathcal A^*\bigl(S_{\bullet}(X_{j})\bigr)=\prod_{j\in J}\negthinspace{}^\downarrow \apl (X_{j}).
\end{multline*}
Since a fibered product of formal CDGA's is clearly formal, we obtain that a wedge of formal spaces is formal, too.

Further examples of formal spaces can be found in geometry. Given a compact, connected Lie group $G$, let $K$ denote the connected component of its neutral element $e$, in the subgroup of elements fixed by a given involution.  The quotient $G/K$, which is a \emph{symmetric space}, is then a formal space, as proved in \cite {cartan}.  Furthermore, Deligne, Griffiths, Morgan and Sullivan showed in \cite {dgms} that compact K\"ahler manifolds are also formal.

It is easy to construct an example of a nonformal CDGA. Let $A=(\Lambda (u,v,w),d)$, where $|u|=|v|=3$ and $|w|=5$ and where $dw=uv$.  Then
\begin{equation*}
\hh _n(A)=\begin{cases} \mathbb Q&:\quad n= 0, 11\\
					   \mathbb Q\oplus \mathbb Q&:\quad n=3,8\\
					   0&:\quad \text {else,}
		\end{cases}
\end{equation*}
where the classes in degree $8$ are represented by $uw$ and $vw$ and the class in degree $11$ by $uvw$.  If $\vp:A\to\hh_*(A)$ is a CDGA map, then $\vp (w)=0$ for degree reasons, which implies that $\vp (uw)=0=\vp (vw)$, since $\vp $ is an algebra map.  Consequently, $\vp$ cannot be a quasi-isomorphism.

\subsection{Models of fiber squares}

The Sullivan model is especially well adapted to studying fibrations.  In particular, as expressed more precisely in the next theorem, ``the Sullivan model of a fiber is the cofiber of the model.''

\begin{theorem}\label{thm:modfib} Let $p:E\to B$ be a Serre fibration such that $B$ is simply connected and $E$ is path connected. Let $F$ denote the fiber of $p$. Suppose that $B$ or $F$ is of finite rational type. 
\begin{enumerate}
\item Given a Sullivan model $\mu :(\Lambda V,d)\to\apl (B)$, let 
$$\xymatrix{ (\Lambda V,d)\ar [rr]^{\apl (p)\circ \mu }\ar [dr]_\iota&&\apl (E)\\
&(\Lambda V\otimes \Lambda W, D)\ar [ur]_{\mu'}^\simeq }$$
be a factorization of $\apl (p)\circ \mu$ as a relative Sullivan algebra, followed by a quasi-isomorphism.  Let $(\Lambda W,\overline D)=\mathbb Q\otimes _{(\Lambda V, d)} (\Lambda V\otimes \Lambda W, D)$, and let $\xymatrix@1{\mu '':(\Lambda W,\overline D)\ar [r]&\apl (F)}$ denote the induced map.  Then $\mu ''$ is a quasi-isomorphism, i.e., there is a commuting diagram in $\cdga$
$$\xymatrix{(\Lambda V,d)\ar[d]^\mu_{\simeq}\ar [r]^(0.35)\iota &(\Lambda V\otimes \Lambda W, D)\ar [d]^{\mu '}_{\simeq}\ar [r]^(0.6)\rho&(\Lambda W, \overline D)\ar [d]^{\mu ''}\\
\apl (B)\ar [r]^{\apl (p)}&\apl (E)\ar [r]^{\apl (j)}&\apl (F)}$$
where $\rho$ is the quotient map and $j$ is the inclusion ma
 
\item Given a Sullivan model  $\xymatrix@1{\mu :(\Lambda V,d)\ar [r]&\apl (B)}$ and a Sullivan minimal model $\xymatrix@1{\mu '':(\Lambda W,d)\ar [r]&\apl (F)}$, there is a relative Sullivan algebra 
$$\iota :(\Lambda V,d)\longrightarrow (\Lambda V\otimes\Lambda W, D)$$  such that $(\Lambda W,d)\cong \mathbb Q\otimes _{(\Lambda V, d)} (\Lambda V\otimes \Lambda W, D)$
and a quasi-isomorphism of cochain algebras 
$$\xymatrix@1{\mu ':(\Lambda V\otimes\Lambda W, D)\ar [r]&\apl (E)}$$
 such that the diagram in (1) commutes, i.e., $E$ has a Sullivan model that is a twisted extension of a Sullivan model of the base by a Sullivan model of the fiber.
 \end{enumerate} 
 \end{theorem}
 
We refer the reader to Proposition 15.5 in \cite{fht} for the proof of the theorem above.  We remark only that the minimality of $(\Lambda W,d)$ in (2) is absolutely essential.

\begin{example} Let  $\xymatrix@1{\Om S^n\ar [r]&PS^n\ar [r]^p&S^n}$ be the based path-space fibration, where $n$ is odd.  Let $\xymatrix@1{\mu : (\Lambda u,0)\ar [r]&\apl (S^n)}$ be the Sullivan model of $S^n$, and consider the relative Sullivan algebra $\xymatrix@1{(\Lambda u,0)\ar [r]&(\Lambda (u,v), d)}$, where $|v|=n-1$ and $dv=u$.  The cochain algebra $(\Lambda (u,v), d)$ is clearly acyclic, as is $\apl (PS^n)$, which implies that $\apl (p)\circ \mu$ extends over $(\Lambda (u,v), d)$ to a quasi-isomorphism of cochain algebras $\xymatrix@1{\mu':(\Lambda (u,v),d)\ar [r]& \apl (PS^n)}$. By Theorem \ref{thm:modfib}, the induced cochain algebra map $\xymatrix@1{\mu'':(\Lambda v,0)\ar [r]&\apl (\Om S^n)}$ is a quasi-isomorphism, which implies that $\hh ^* (\Om S^n;\mathbb Q)\cong \mathbb Q[v]$, when $n$ is odd.

More generally, consider the based path-space fibration $\xymatrix@1{\Om X\ar [r]&PX\ar [r]^p&X}$, where $X$ is a simply connected space.  Suppose that $\xymatrix@1{\mu :(\Lambda V, d)\ar [r]&\apl (X)}$ is a Sullivan model of $X$.  

Using notation that is standard in rational homotopy theory, let $\overline V$ be the graded $\mathbb Q$-vector space that is the suspension of $V$, i.e., $\overline V^n=V^{n+1}$.  Let $S$ be the derivation of $\Lambda (V\oplus \overline V)$ specified by $S(v)=\bar v$ and $S(\bar v)=0$ for all $v\in V$. Define $\big(\Lambda (V\oplus \overline V),D\big)$ by $D\bar v= -S(dv)$, which implies that $D\bar v\in \Lambda V\otimes \Lambda ^+\overline V$, i.e., each summand of $D\bar v$ contains at least one factor from $\overline V$.  The CDGA $\bigl(\Lambda (V\oplus \overline V),D\bigr)$  is easily seen to be acyclic.

We can now define a morphism of CDGA's 
$$\mu ':\bigl(\Lambda (V\oplus \overline V),D\bigr)\to\apl (PX)$$ 
by setting $\mu'(v)=\mu (v)$ for all $v\in V$ and $\mu'(\bar v)=0$ for all $\bar v\in \overline V$ and extending multiplicatively.  Since both the source and the target of $\mu$ are acyclic, $\mu'$ is a quasi-isomorphism.  We have therefore constructed a commutative diagram of CDGA's
$$\xymatrix{
(\Lambda V,d)\ar [d]_{\mu}^\simeq \ar[r]^{\iota}&\bigl(\Lambda (V\oplus \overline V),D\bigr)\ar[d]_{\mu'}^\simeq\\
\apl(X)\ar [r]^p&\apl(PX)}$$
to which we can apply Theorem \ref{thm:modfib}. The induced morphism of CDGA's 
$$\xymatrix@1{\mu '':(\Lambda \overline V, 0)\ar [r]&\apl (\Om X)}$$ 
is  therefore a quasi-isomorphism, so that $\hh ^*(\Om X;\mathbb Q)\cong \Lambda \overline V$. 
\end{example}

Theorem \ref{thm:modfib} is a consequence of the following, more general result concerning fiber squares, for which the slogan is ``the model of the pullback is the pushout of the models.''

\begin{theorem} \label{thm:modpb}  Let $p:E\to B$ be a Serre fibration, where $E$ is path connected and $B$ is simply connected, with fiber $F$.  Let $f:X\to B$ be a continuous map, where $X$ is simply connected. Suppose that $B$ or $F$ is of finite rational type.  Consider the pullback
$$\xymatrix{E\underset B \times X\ar [r]^{\bar f}\ar [d]^{\bar p}&E\ar [d]^p\\
X\ar [r]^f&B.}$$
Given a commuting diagram of CDGA's
$$\xymatrix{(\Lambda U,d)\ar [d]_\simeq ^\nu&(\Lambda V,d)\ar [l]_\vp\ar [r]^\iota \ar [d]_\simeq ^\mu &(\Lambda V\otimes \Lambda W, D)\ar [d]_\simeq ^{\mu '}\\
\apl (X)&\apl (B)\ar [l]_{\apl(f)}\ar [r]^{\apl (p)}&\apl (E),}$$
consider the pushout diagram in $\cdga$
$$\xymatrix{
(\Lambda V,d)\ar[d]^\iota \ar[r]^{\varphi}&(\Lambda U, d)\ar[d]\\
(\Lambda V\otimes \Lambda W, D)\ar [r]&(\Lambda U\otimes \Lambda W, \bar D).}$$
Then the induced map of cochain algebras
$$\xymatrix{(\Lambda U\otimes \Lambda W, \bar D)\ar [r]& \apl (E\underset B\times X)}$$
is a quasi-isomorphism of CDGA's.
\end{theorem}

We again refer the reader to \cite {fht} for the proof of this theorem, in the guise of their Proposition 15.8.

\begin{example}  Let $X$ be a simply connected space of finite rational type, with Sullivan model $\xymatrix@1{\varphi:(\Lambda V,d)\ar [r]& \apl (X)}$.  Consider the free-loop pullback square
$$\xymatrix{\mathcal L X\ar [d]^e\ar [r] ^j&X^I\ar [d] ^{p}\\ 
		   X\ar[r]^(0.40)\Delta &X\times X,}$$
where $\Delta$ is the diagonal map and where $p(\lambda )=\big(\lambda (0),\lambda (1)\big)$ for all paths $\lambda:I\to X$.   The free loop space $\mathcal LX$ is thus $\{\lambda \in X^I\mid \lambda (0)=\lambda (1)\}$, which is homeomorphic to $X^{S^1}$, the space of unbased loops on $X$.

It is easy to check that 
$$\xymatrix{(\Lambda V,d)\otimes (\Lambda V,d)\ar [d]^{\widehat\varphi}_{\simeq} \ar [r]^(0.6)m&(\Lambda V,d)\ar [d]^\varphi_{\simeq}\\
			\apl (X\times X)\ar [r]^(0.55){\apl(\Delta)}&\apl (X)}$$
commutes, where $m$ is the multiplication map on $(\Lambda V, d)$ and $\widehat\varphi$ is the composite
$$(\Lambda V,d)\otimes (\Lambda V,d)\xrightarrow {\varphi\otimes \varphi}\apl(X)\otimes \apl (X)\xrightarrow {\alpha _{X,X}}\apl (X\times X).$$

Furthermore
$$\xymatrix{(\Lambda V,d)\otimes (\Lambda V,d)\ar [d]^{\widehat\varphi}_{\simeq} \ar [r]^(0.45)\iota&(\Lambda (V'\oplus V''\oplus \overline V),D)\ar [d]^{\Phi}_{\simeq}\\
			\apl (X\times X)\ar [r]^{\apl(p)}&\apl (X^{I})}$$
commutes as well, where $\iota $ is a relative Sullivan algebra ($V'$ and $V''$ are two copies of $V$), $\Phi$ is an appropriate extension of $\widehat \varphi$ and $D$ is specified as follows. 

Let $S: \Lambda (V'\oplus V''\oplus \overline V)\longrightarrow  \Lambda (V'\oplus V''\oplus \overline V)$ be the derivation of degree $-1$ specified by $S(v')=\bar v=S(v')$ and $S(\bar v)=0$.  Then 
$$D(\bar v):=v''-v'-\sum _{n\geq 1} \frac {(SD)^n}{n!} (v').$$

Applying Theorem \ref{thm:modpb} we obtain as a Sullivan model of $\mathcal LX$
$$\bigl (\Lambda (V\oplus \overline V), \overline D)\cong (\Lambda V,d )\underset {(\Lambda V, d)^{\otimes 2}}\otimes(\Lambda (V'\oplus V''\oplus \overline V),D),$$
i.e., $\overline D (\bar v)= -\overline S(dv)$, where  $\overline S(v)=\bar v$ and $\overline S(\bar v)=0$.

Sullivan and Vigu\'e used this model to prove that if $\hh ^*(X;\mathbb Q)$ requires at least two algebra generators, then the rational Betti numbers of $\hh^*(\mathcal LX;\mathbb Q)$ grow exponentially, which implies in turn that $X$ admits an infinite number of distinct closed geodesics, when $X$ is a closed Riemannian manifold \cite{suvp}.

It is well known (cf., e.g., \cite{loday}) that  $\hh ^*(\mathcal LX;\mathbb Z)\cong \operatorname{HH}_{*}\big(C^*(X;\mathbb Z)\big)$, the Hoch\-schild homology of the cochains on $X$.  The Sullivan model of the free loop space constructed above is thus also a tool for understanding Hochschild homology.
\end{example}

\subsection {Lusternik-Schnirelmann category}

One of the most spectacular successes of the Sullivan minimal model
has been in its application to studying and exploiting the numerical homotopy
invariant known as Lusternik-Schnirelmann (L.-S.) category.  

\begin{definition} A \emph{categorical covering} of a space $X$ is an open cover of $X$
such that each member of the cover is contractible in $X$. The
\emph{L.-S. category} of a topological space $X$, denoted $\cat X$, is equal to
$n$ if the cardinality of the smallest categorical covering of $X$ is
$n+1$. \end{definition}

L.-S. category is in general extremely difficult to compute.  It is
trivial, however, to prove that the L.-S. category of a contractible space is 0 and that $\cat S^n=1$ for all $n$. Similarly, the L.-S. category of any suspension is $1$.  More generally, a space $X$ is a co-H-space if and only if $\cat X\leq 1$.  

The proof of this last equivalence is most easily formulated in terms of an equivalent definition of L.-S. category, which requires the following construction, due to Ganea.

\begin{definition}  Let $p:E\longrightarrow  X$ be a fibration over a based topological space $(X,x_{0})$.  Let $j:F\hookrightarrow E$ denote the inclusion of the fiber of $p$ over $x_{0}$, with mapping cone $C_{j}$.  Let $\hat p: C_{j}\longrightarrow  X$ denote the induced continuous map, which can be factored naturally as a homotopy equivalence followed by a fibration:
$$\xymatrix{C_{j}\we&E' \fib^{p'}&X.}$$
There is then a commutative diagram
$$\xymatrix{E\ar [rr]\ar[dr]_{p}&&E'\ar[dl]^{p'}\\ &X}$$
called  the \emph{fiber-cofiber construction} on $p$.

Let $p:PX\longrightarrow  X$ denote the (based) path fibration over $X$.  Iterating the fiber-cofiber construction repeatedly leads to a commutative diagram
$$\xymatrix{P_{0}X=PX\ar[d]_{p}\ar [r]&P_{1}X\ar [dl]^(0.35){p_{1}}\ar [r]&P_{2}X\ar [dll]^(0.4){p_{2}}\ar [r]&\cdots\ar [r]&P_{n}X \ar [dllll]^{p_{n}}\ar [r]&\cdots     \\ X,    }$$
in which $P_{n}X$ is the $n^{th}$ \emph{Ganea space} for $X$ and $p_{n}$ is the $n^{th}$ \emph{Ganea fibration}.
\end{definition}

\begin{proposition} If $X$ is a normal space, then $\cat X\leq n$ if and only if the fibration $p_{n}$ admits a section.
\end{proposition}

For a proof of this proposition, we refer the reader to Proposition 27.8 in \cite{fht}.

Another equivalent definition of the L.-S. category of a based space $(X,x_0)$ is expressed in terms of the \emph{fat wedge} on $X$ 
$$T^nX:=\{(x_1,...,x_n)\in X^n\mid \exists\; i \text{ such that } x_i=x_0\}.$$

\begin{proposition}\label{prop:whcat}  If $X$ is a path-connected CW-complex, then the following conditions are equivalent.
\begin{enumerate}
\item $\cat X\leq n$.
\item The iterated diagonal $\xymatrix@1{\Delta ^{(n)}: X\ar [r]& X^{n+1}}$ factors up to homotopy through the fat wedge $T^{(n+1)}X$, i.e., there is a map $\xymatrix@1{\delta :X\ar[r]& T^{n+1}X}$ such that  the diagram
$$\xymatrix{X\ar [dr]^\delta \ar [rr]^{\Delta ^{(n)}}&& X^{n+1}\\
                                 &T^{n+1}X\ar [ur]^i}$$
commutes up to homotopy.
\end{enumerate}
\end{proposition}

We refer the reader to Proposition 27.4 in \cite {fht} for the proof of this equivalence.

Applying Proposition \ref{prop:whcat}, we obtain the following useful upper bound on L.-S. category (cf., Proposition 27.5 in \cite {fht}).

\begin{corollary}  If $X$ is an $(r-1)$-connected CW-complex of dimension $d$, where $r\geq 1$, then $\cat X\leq d/r$.
\end{corollary}

On the other hand, a lower bound on $\cat X$ is given by the \emph{cuplength} $c(X)$ of  $H^*(X;\mathbb Q)$, i.e, the greatest integer $n$ such that there exist $a_{1},...,a_{n}\in H^*(X;\mathbb Q)$ satisfying $a_{1}\cup \cdots \cup a_{n}\neq 0$.  We leave it as an easy exercise to prove that $c(X)\leq \cat X$ for all path-connected, normal spaces $X$.

\begin{example} Observe that $c(\mathbb CP^n)=n$, so that $\cat \mathbb CP^n\geq n$. On the other hand, $\mathbb CP^n$ is 1-connected and of dimension $2n$, implying that $\cat \mathbb CP^n\leq 2n/2=n$.  Thus, $\cat \mathbb CP^n=n$.\end{example}\medskip

Within the realm of rational homotopy theory, it makes sense to consider the following invariant derived from L.-S. category.

\begin{definition}  The \emph{rational category} of a simply connected space $X$, denoted $\cat_{0}X$, is defined by
 $$\cat_{0}X:=\min \{\cat Y\mid \text{ $X$ and $Y$ have the same rational homotopy type}\}.$$\end{definition}
 
 As proved in \cite {fht}  (Proposition 28.1), if $X$ is a simply connected CW-complex, then $\cat _{0}X=\cat X_{0}.$ Furthermore, it is obvious that $\cat _{0}X\leq \cat X$ for all $X$.  As we show below, this inequality is sharp, i.e., there are spaces $X$ for which $\cat _{0}X= \cat X$.  On the other hand, the inequality can certainly be strict, as the case of a mod $p$ Moore space easily illustrates.

The next theorem has turned out to be a crucial tool in proving numerous significant results in rational homotopy theory, such as many of the dichotomy theorems (cf., Section 2.4).

\begin{theorem}[The Mapping Theorem]\label{thm:mapthm} Let $f:X\longrightarrow  Y$ be a continuous map between simply connected spaces.  If $\pi _{*}f\otimes \mathbb Q$ is injective, then $\cat _{0}X\leq \cat _{0}Y$.
\end{theorem}

The original proof of the Mapping Theorem relied on Sullivan models.  There is now a purely topological and relatively simple proof, which is given in \cite {fht} (Theorem 28.6).

As a first application of the Mapping Theorem, we mention the amusing and useful corollary below, which follows immediately from the fact that the natural map from the $(n+1)^{st}$ Postnikov fiber to the $n^{th}$ Postnikov fiber of a space induces an injection on homotopy groups.

\begin{corollary} Let $X$ be a connected CW-complex.  Let $X(n)$ denote the $n^{th}$ Postnikov fiber of $X$, for all $n\geq 1$. Then
$$\cdots \leq\cat _{0}X(n+1)\leq \cat _{0}X(n)\leq\cdots \leq \cat_{0}X(2)\leq \cat _{0}X.$$
\end{corollary}

One great advantage of rational category, as opposed to the usual L.-S. category, is that it can explicitly calculated in terms of the Sullivan model, as stated in the next theorem.  

\begin{theorem}\label{thm:ratcat} Let $\xymatrix@1{\varphi:(\Lambda
V,d)\we &\apl (X)}$ be the Sullivan minimal model of a simply connected space $X$ of finite rational type.  Let $(\Lambda V/\Lambda
^{>n}V,\bar d)$ denote the CDGA obtained by taking the quotient of
$(\Lambda V,d)$ by the ideal of words of length greater than $n$, and
let 
\begin{equation}
\xymatrix{ (\Lambda V,d)\ar[rr]^{q}\ar[dr]_{i}&&(\Lambda
V/\Lambda ^{>n}V,\bar d)\\
&(\Lambda (V\oplus W),d)\ar [ur]_{p}^{\simeq}}
\end{equation}
be a factorization of the quotient map $q$ as a relative Sullivan algebra, followed by a surjective quasi-isomorphism.
Then $cat_{0} X\leq n$ if and only if $i$ admits a CDGA retraction
$\xymatrix@1{\rho : (\Lambda (V\oplus W),d)\ar [r]&(\Lambda V,d)},$
i.e.,  $\rho i=Id_{(\Lambda V,d)}$.
\end{theorem}

The fat wedge formulation of the definition of L.-S. category is crucial in the proof of this theorem, for which we refer the reader to Propositions 29.3 and 29.4 in \cite{fht}.

\begin{example} Since $H^n(S^n;\mathbb Q)=\mathbb Q$, the rationalization $S^n_{0}$ of the $n$-sphere is not contractible and therefore $\cat _{0}S^n=\cat S^n_{0}>0$.  On the other hand, $\cat _{0}S^n\leq \cat S^n=1$, whence $\cat _{0}S^n=1$, providing the promised example of equality between rational category and L.-S. category of a space. 

This calculation can also be carried out easily using the Sullivan model (cf., Section 2.1). If $n$ is odd, then the Sullivan model is $(\Lambda (x), 0)$, where $x$ is of degree $n$. Observe that  $\Lambda (x)/\Lambda ^{>1}(x)$ is isomorphic to $\Lambda (x)$, since $x$ is of odd degree. Since the quotient map $q$ is itself the identity map in this case, it follows trivially that $\cat _{0}S^n\leq 1$.

If $n$ is even, the relevant Sullivan model is $(\Lambda (y,z),d)$, where $\deg y=n$, $\deg z=2n-1$ and $dz=y^2$.  An easy calculation shows that 
$$(\Lambda (y,z)/\Lambda ^{>1} (y,z),\bar d)=(\mathbb Q\oplus \mathbb Q\cdot y\oplus \mathbb Q\cdot z, 0).$$
It is not too difficult to show that the quotient map $q$ factors as
$$\xymatrix{ (\Lambda (y,z),d)\ar[rr]^{q}\ar[dr]_{i}&&(\mathbb Q\oplus \mathbb Q\cdot y\oplus \mathbb Q\cdot z,,0)\\
&(\Lambda (y,z)\otimes \Lambda W,D)\ar [ur]_{p}^{\simeq}}$$
where $DW\subset \Lambda (y,z)\otimes \Lambda ^+W$, i.e., the differential of any generator of $W$, if nonzero, is a sum of words, all of which contain at least one letter from $W$.  We can therefore define a CDGA retraction $\rho$ by setting $\rho (w)=0$ for all $w\in W$, implying that $\cat _{0}S^n\leq 1$.
\end{example}

Though Theorem \ref{thm:ratcat} does simplify the calculation of rational category by making it purely algebraic, the computations involved are still difficult, which led Halperin and Lemaire to propose the following, apparently weaker numerical invariant of rational homotopy \cite {hallem}.

\begin{definition}   Let $X$ be a simply connected space of finite rational type, with Sullivan model $\xymatrix@1{\varphi:(\Lambda
V,d)\we &\apl (X)}$.  If the map $i$ in diagram (2.1) of Theorem \ref{thm:ratcat} admits a retraction as morphisms in the category of $(\Lambda V,d)$-modules, then $\mcat _{0}X\leq n$.
\end{definition}

As it turned out, however, the apparent weakness of $\mcat_{0}$ was only an illusion.

\begin {theorem}\label{thm:mine} $\mcat_{0}X=\cat _{0}X$ for all
simply connected spaces $X$ of finite
rational type.
\end{theorem}

For the proof of this theorem, which requires a deep understanding of the factorization of the quotient map, we refer the reader to Theorem 29.9 in \cite {fht}.

Theorem \ref{thm:mine} implies that to show that $\cat _{0}X\leq n$, it suffices to find a $(\Lambda V, d)$-module retraction of $i$ in diagram (2.1), which has proven to be a very effective simplification. We next outline briefly one application of this simplification, to the study of the additivity of L.-S. category.

It is not difficult to show that $\cat (X\times Y)\leq
\cat  X+\cat Y$ for all normal spaces $X$ and $Y$.  At the end of the 1960's Ganea observed that
in the only known examples for which $\cat (X\times Y)\not=\cat X+\cat
Y$, the spaces $X$ and $Y$ had homology torsion at distinct primes. 
He conjectured therefore that $\cat (X\times S^n)=\cat X+ 1$ for all
spaces $X$ and all $n\geq 1$, since $S^n$ has no homology torsion
whatsoever.

In fact, as stated precisely below, if we forget torsion completely and work rationally, then L.-S. category is indeed additive.

\begin{theorem} If $X$ and $Y$ are simply connected topological spaces of finite rational type, then $\cat_{0}(X\times Y)=\cat _{0}X+\cat _{0}Y.$
\end{theorem}

The proof of this theorem depends in an essential way on Theorem \ref{thm:mine}.  We refer the reader to Sections 29(h) and 30(a) in \cite{fht}Ê for further details.

As an epilogue to this story of Sullivan minimal models and L.-S.
category, we mention that in 1997 Iwase applied classical
homotopy-theoretic methods to the construction of a counter-example to
Ganea's conjecture \cite {iwase}.  In particular, he built a $2$-cell complex $X$ such that $\cat (X\times S^n)=2=\cat X$.

\subsection {Dichotomy}

There is a beautiful dichotomy governing finite CW-com\-plex\-es in rational homotopy theory, expressed as follows: the rational homotopy groups of a finite CW-complex  either are of finite total dimension as graded rational vector space or grow exponentially.  We first examine the former case, that of \emph{elliptic spaces}, then the latter case, that of \emph{hyperbolic spaces}.

\begin{definition} A simply connected topological space $X$ is \emph{rationally elliptic} if 
$$\dim \hh ^*(X;\mathbb Q)<\infty\text{ and }\dim \pi_*(X)\otimes \mathbb Q<\infty.$$  The \emph{formal dimension} $\fdim X$ of a rationally elliptic space $X$ is defined by
$$\fdim X:=\max\{k\mid \hh ^k(\Lambda V,d)\not= 0\}.$$  
The \emph{even exponents} of a rationally elliptic space $X$ are the postive integers $a_1,...,a_q$ such that there is a basis $(y_j)_{1\leq j\leq q}$ of $\pi _{\text{even}}X\otimes \mathbb Q$ with $\deg y_j=2a_j$.  Similarly, the \emph{odd exponents} of $X$ are the positive integers $b_1,...,b_p$ such that there is a basis $(x_i)_{1\leq i\leq p}$ of $\pi _{\text{odd}}X\otimes \mathbb Q$ with $\deg x_i=2b_i-1$.

A minimal Sullivan algebra $(\Lambda V,d)$ is \emph{elliptic} if the associated rational space $<(\Lambda V,d)>$ is elliptic.  The \emph{formal dimension} of an elliptic minimal Sullivan algebra is the formal dimension of the associated space.  
\end{definition}

\begin{example}  Spheres, complex projective spaces, products of elliptic spaces, and homogeneous spaces are examples of elliptic spaces.
\end{example}

 The following special case of elliptic spaces is important for understanding general elliptic spaces.
 
 \begin{definition} A minimal Sullivan algebra $(\Lambda V,d)$ is \emph{pure} if $\dim V<\infty$, $d|_{V^{\text{even}}}=0$ and $d(V^{\text{odd}})\subseteq \Lambda V^{\text{even}}$. A space is \emph{pure} if its Sullivan model is pure. 
 \end{definition}
 
Note that a pure space $X$ is elliptic if and only if  $\dim \hh ^*(X;\mathbb Q)<\infty$.
 
A pure Sullivan algebra $(\Lambda V, d)$ admits a differential filtration 
$$F_k(\Lambda V,d)=\Lambda V^{\text{even}}\otimes \Lambda ^{\leq k}V^{\text{odd}}.$$  
In particular, $dF_k(\Lambda V,d)\subset F_{k-1}(\Lambda V,d)$.  Write
$$\hh_k(\Lambda V,d)=\frac {\ker \big(d:F_k(\Lambda V,d)\longrightarrow F_{k-1}(\Lambda V,d)\big)}{\im \big(d:F_{k+1}(\Lambda V,d)\longrightarrow F_k(\Lambda V,d)\big)}.$$

The following list of the most important properties of pure Sullivan algebras summarizes Propositions 32.1 and 32.2 in \cite {fht}.

\begin{proposition}  Let $(\Lambda V,d)$ be a pure, minimal Sullivan algebra..
\begin{enumerate}
\item $\dim \hh ^*(\Lambda V,d)<\infty \Leftrightarrow \dim H_0(\Lambda V,d)<\infty$.
\smallskip
\item If $\dim  \hh ^*(\Lambda V,d)<\infty $, then $\hh^n(\Lambda V,d)$ is a $1$-dimensional subspace of $\hh_r(\Lambda V,d)$, where $n$ is the formal dimension of $(\Lambda V, d)$ and $r=\max\{ k\mid \hh _k(\Lambda V,d)\not =0\}$.
\smallskip
\item If $\dim  \hh ^*(\Lambda V,d)<\infty $, then $r=\dim V^{\text{odd}}-\dim V^{\text{even}}$.  Thus,  
$$\dim V^{\text{odd}}=\dim V^{\text{even}} \Leftrightarrow \hh^*(\Lambda V,d)=\hh _0(\Lambda V,d).$$
\smallskip
\item If $\dim  \hh ^*(\Lambda V,d)<\infty $ and $(\Lambda V, d)$ is a Sullivan minimal model of $X$, then 
$$\fdim X=\sum _i(2b_i-1)-\sum _j (2a_j-1),$$ 
where the even and odd exponents of $X$ are $a_{1},..,a_{q}$ and $b_{1},...,b_{p}$, respectively.
\end{enumerate}
\end{proposition}

We present next a tool for determining whether spaces are elliptic, based on the notion of pure spaces.

\begin{definition} Let $(\Lambda V,d)$ be a minimal Sullivan algebra such that $\dim V<\infty$. Filter $(\Lambda V,d)$ by 
$$\mathcal F^p(\Lambda V,d)=\bigoplus_ {k+l\geq p} (\Lambda V^{\text{even}}\otimes \Lambda ^kV^{\text{odd}})^l.$$
The induced spectral sequence is called the \emph{odd spectral sequence} and converges to $\hh^*(\Lambda V,d)$.
\end{definition}

Observe that the $E_0$-term of the odd spectral sequence is the associated graded of $(\Lambda V, d_\sigma)$, where $d_\sigma (V^{\text{even}})=0$, $d_\sigma (V^{\text{odd}})\subseteq \Lambda V^{\text{even}}$ and $(d-d_\sigma)(V^{\text{odd}})\subseteq \Lambda V^{\text{even}}\otimes \Lambda ^+ V^{\text{odd}}$.  We call $(\Lambda V, d_\sigma)$ the \emph{associated pure Sullivan algebra} of $(\Lambda V,d)$.

\begin{proposition} Under the hypotheses of the definition above,
$$\dim \hh^*(\Lambda V,d)<\infty\Leftrightarrow \dim \hh^*(\Lambda V,d_\sigma )<\infty.$$
\end{proposition}

Thus, $(\Lambda V,d)$ is elliptic if and only if $(\Lambda V,d_\sigma)$ is elliptic.

\begin{proof} As the odd spectral sequence converges from $ \hh^*(\Lambda V,d_\sigma )$ to $ \hh^*(\Lambda V,d)$, one implication  is clear.  An algebraic version of the Mapping Theorem (Theorem \ref{thm:mapthm}) plays an essential role in the rest of the proof.  We refer the reader to Proposition 32.4 in \cite {fht} for the complete proof.
\end{proof}

\begin{example} Consider $(\Lambda V,d)=(\Lambda (a_2, x_3, u_3, b_4, v_5, w_7),d)$, where the subscript of a generator equals its degree and $d$ is specified by $da=0$, $dx=0$, $du=a^2$, $db=ax$, $dv=ab-ux$ and $dw=b^2-vx$.  Its associated pure Sullivan algebra is $(\Lambda (a, x, u, b, v, w),d_\sigma)$, where $d_\sigma a=0$, $d_\sigma x=0$, $d_\sigma u=a^2$, $d_\sigma b=0$, $d_\sigma v=ab$ and $d_\sigma w=b^2$.  A straightforward calculation shows that 
$$\hh^*(\Lambda V,d_\sigma)=\mathbb Q\cdot a\oplus \mathbb Q\cdot b \oplus \Lambda y/(y^3)  \oplus \Lambda z/(z^3),$$
where $y$ is represented by $bu-av$ and $z$ is represented by $aw-bv$.  In particular $\dim \hh^*(\Lambda V, d_\sigma)<\infty$, which implies that $(\Lambda V,d)$ is elliptic.
\end{example}

The next theorem describes the amazing numerology of elliptic spaces and Sullivan algebras, which imposes formidable constraints on their form.  

\begin{theorem} Let $(\Lambda V,d)$ be an elliptic, minimal Sullivan algebra of formal dimension $n$ and with even and odd exponents $a_1,...,a_q;b_1,...,b_p$.  Then:
\begin{enumerate}
\item $\sum _{i=1}^p (2b_i-1)-\sum _{j=1}^q (2a_j-1)=n$;
\smallskip
\item $\sum _{j=1}^q 2a_j\leq n$;
\smallskip
\item $\sum _{i=1}^p(2b_1-1)\leq 2n-1$; and
\smallskip
\item $\dim V^{\text{even}}\leq \dim V^{\text{odd}}$.
\end{enumerate}
\end{theorem}

As a consequence of this theorem, we know, for example, that if $(\Lambda V,d)$ is an elliptic, minimal Sullivan algebra of formal dimension $n$, then $V=V^{\leq 2n-1}$, $\dim V^{>n}\leq 1$ and $\dim V\leq n$.

\begin{proof} One first proves by induction on $\dim V$ that the formal dimensions of $(\Lambda V,d)$ and of its associated pure Sullivan algebra are the same, reducing the proof of the theorem to the pure case. For further details, we refer the reader to Theorem 32.6 in \cite {fht}.
\end{proof}

\begin{definition} The \emph{Euler-Poincar\'e characteristic} of a graded vector space $W$ of finite type is the integer
$$\chi _W=\sum _i(-1)^{i}\dim W^i=\dim W^{\text{even}}-\dim W^{\text {odd}}.$$
\end{definition}
 
 It is easy to show that $\chi _W=\chi _{\hh^*(W,d)}$, for any choice of differential $d$ on $W$.
 
 \begin{proposition}\label{prop:eulchar}  If $(\Lambda V,d)$ is an elliptic, minimal Sullivan algebra, then $\chi_{\Lambda V}\geq 0$ and $\chi _V\leq 0$.  Furthermore, the following statements are equivalent.
 \begin{enumerate}
 \item $\chi _{\Lambda V}>0$.
 \smallskip
 \item $\hh ^*(\Lambda V, d)=\hh ^{\text{even}} (\Lambda V,d)$.
 \smallskip
 \item $\hh ^*(\Lambda V,d)=\Lambda (y_1,...,y_q)/(u_1,..., u_p),$ where $(u_1,...,u_p)$ is a regular sequence.
 \smallskip
 \item $(\Lambda V,d)$ is isomorphic to a pure, minimal Sullivan algebra.
 \smallskip
 \item $\chi _V=0$.
 \end{enumerate}
 \end{proposition}

 \begin{proof} The proof of this proposition reduces essentially to Poincar\'e series calculations.  We refer the reader to Proposition 32.10 in \cite {fht} for details of the calculations.
 \end{proof}
 
 \begin{definition} Let $X$ be an elliptic space with minimal Sullivan model $(\Lambda V,d)$.  The \emph{homotopy Euler characteristic} $\chi_{\pi}(X)$ of $X$ is defined to be $\chi_{V}$.
 \end{definition}
 
 Propostion \ref{prop:eulchar} implies that $\chi_{\pi}(X)\leq 0$ for all elliptic spaces $X$ and that $\chi_{\pi}(X)=0$ if and only if $X$ is a pure elliptic space.

\begin{example}[Application to free torus actions (Example 3 in section 32(e) of \cite {fht})] Let $T$ denote the $r$-torus, i.e., the product of $r$ copies of $S^1$.  Suppose that $T$ acts smoothly and freely on a simply connected, compact, smooth manifold $M$.  There exists then a smooth principal bundle $M\longrightarrow M/T$ and thus a classifying map $M/T\longrightarrow BT$ with homotopy fiber $M$.

If $M$ is elliptic, then $M/T$ is also elliptic, since $M/T$ is compact and $BT=(\mathbb CP^\infty)^r$.  Furthermore,
$$0\geq \chi _\pi (M/T)=\chi _\pi (M)+\chi _\pi (BT)=\chi _\pi (M)+r ,$$
implying that $r\leq - \chi _\pi (M)$.
\end{example}

Now we go to the other extreme.

\begin{definition} A simply connected space $X$ with the homotopy type of a finite CW-complex is \emph{rationally hyperbolic} if $\dim \pi _*(X)\otimes \mathbb Q=\infty$.
\end{definition}

The following theorem, which justifies the terminology ``hyperbolic,'' is Theorem 33.2 in \cite {fht}.  Its proof depends strongly on the Mapping Theorem (Theorem \ref{thm:mapthm}).

\begin{theorem}  If $X$ is a rationally hyperbolic space, then there exist $C>1$ and $N\in \mathbb Z$ such that 
$$\sum _{i=0}^n\dim \pi _i(X)\otimes \mathbb Q\geq C^n$$
for all $n\geq N$.  
\end{theorem}

In other words, the rational homotopy groups of $X$ \emph{grow exponentially}.  Moreover, as stated more precisely in the next theorem (Theorem 33.3 in \cite {fht}), there are no ``long gaps" in the rational homotopy groups of $X$.

\begin{theorem}  If $X$ is a rationally hyperbolic space of formal dimension $n$, then for all $k\geq 1$, there exists $i\in (k,k+n)$ such that $\pi _iX\otimes \mathbb Q\not=0$.  Furthermore, for $k\gg 0$,
$$\sum _{i=k+1}^{k+n-1} \dim \pi _i(X)\otimes \mathbb Q\geq \frac {\dim \pi _kX\otimes \mathbb Q}{\dim \hh^*(X;\mathbb Q)}.$$
\end{theorem}

Consequently, if $X$ has formal dimension $n$, then 
$$\text{$X$ is rationally elliptic}\Leftrightarrow \pi _j(X)\otimes \mathbb Q=0 \quad\forall j\in [2n, 3n-2],$$
a simple and lovely test of ellipticity.

\section {Commutative algebra and rational homotopy theory}

In the late 1970's  two algebraists, Luchezar Avramov, then at the
University of Sofia, and Jan-Erik Roos of the University of Stockholm,
discovered and began to exploit a  deep   connection between local
algebra and rational homotopy theory.  In 1981 they established contact
with the rational homotopy theorists, initiating a powerful synergy that led
to a multitude of important results in both fields.  In this section, of a more expository nature than the preceding sections, we describe certain of the most important results of this collaboration.  For further details we refer the reader to Section 4 of  \cite{hess}.

Roos' interest in rational homotopy theory was inspired by work of Jean-Michel Lemaire
on Serre's question concerning rationality of Poincar\'e series of the rational
homology of loop spaces and by the work of local algebraists on the
analogous question of Kaplansky and Serre for local rings.  More precisely,
Lemaire had studied the Poincar\'e series
$$\sum \sb {n\geq 0} \dim\sb {\mathbb Q} H\sb n(\Omega E;\mathbb Q) \cdot z\sp
n$$ for $E$ a finite, simply connected CW-complex, while local algebraists
were interested in the series
$$\sum \sb {n\geq 0} \dim\sb {\mathbb F} \operatorname{Ext}\sb R\sp n(\mathbb F,
\mathbb F)\cdot z\sp n$$
for $R$ a local ring with residue field $\mathbb F$.  (Here and throughout this section, all local rings are assumed to be commutative and noetherian.)
In both cases, the goal was to determine whether the series always
represented a rational function.

Roos established a
research program to study the homological properties of local rings,
in particular those whose maximal ideal $\mathfrak m$ satisfied $\mathfrak m\sp
3=0$, the first nontrivial case for Poincar\'e series calculations.   He
realized that in order to study local rings, it was useful, or even
necessary, to work in the larger category of (co)chain algebras.   
By 1976
he had proved the equivalence of Serre's problem for  CW-complexes $E$ such
that $\dim E=4$ and of the Kaplansky-Serre problem for local rings $(R,\mathfrak 
m)$ such that $\mathfrak m\sp 3=0$ \cite {roos}.  

Avramov ascribes his original interest in rational homotopy theory to Levin's result from 1965  that if $R$ is a local  ring
with residue field $\mathbb F$, then $\operatorname {Tor}\sp R(\mathbb F,\mathbb F)$
is a graded, divided powers Hopf algebra \cite {levin}.  Since the dual of a graded,
divided powers Hopf algebra is the universal enveloping algebra of a uniquely
defined graded Lie algebra ($\operatorname{char}\mathbb F =0$ \cite {mimo},
$\operatorname{char}\mathbb F >2$
\cite {andrŽ}, $\operatorname{char}\mathbb F =2$ \cite {sjo}), to any local ring $R$ is associated a graded Lie algebra $\pi \sp *(R)$, the \emph{homotopy Lie algebra} of $R$.  
Based on results in characteristic 0 due to Gulliksen in the late
1960's, Avramov proved that if
$R\longrightarrow S$ is a homomorphism of local rings such that $S$ is $R$-flat, then there is an
exact sequence of groups
{\smaller{$$...\to\pi^n R\longrightarrow \pi\sp n(S\otimes \sb R \mathbb F)\longrightarrow \pi\sp n(S)\longrightarrow \pi\sp
n(R)\xrightarrow{\delta ^n}\pi\sp {n+1}(S\otimes \sb R \mathbb F) \longrightarrow\pi ^{n+1}S\to ... .$$}}
The existence of such a long exact sequence of homotopy groups confirmed Av\-ra\-mov's intuition that rational homotopy invariants
were analogous to homology invariants of local rings in
arbitrary characteristic.

In 1980  David Anick
constructed a finite, simply connected CW-complex $E$ of dimension $\leq
4$ such that the Poincar\'e series of the homology of $\Omega E$ was not
rational  \cite {anick1}, \cite {anick2}.  Anick's construction
interested rational homotopy theorists because of its relation to the
dichotomy between elliptic and hyperbolic spaces; see Section 2.4.  Local
algebraists were interested because of Roos's result, which allowed the
transcription of Anick's space into a local ring $(R,\mathfrak m)$ with $\mathfrak
m\sp 3=0$ and with irrational Poincar\'e series. Shortly after Anick's result became known, Roos and his student Clas
L\"ofwall discovered other examples of local rings with irrational
Poincar\'e series that they obtained by completely different methods \cite
{loroos}.

The converging interests of rational homotopy theorists and of local
algebraists led to direct contact between the two groups in 1981. 
Inspired by the work of Roos and his colleagues,  F\'elix and Thomas began to work on
calculating the radius of convergence of the Poincar\'e series of a loop
space \cite {fŽthom},  establishing the following beautiful
characterization:  a
simply connected space 
$E$ of finite L.~-S.~category is rationally elliptic if and only if the radius of convergence of the
Poincar\'e series of $\Omega E$ is 1.  If $E$ is rationally hyperbolic,
then the radius of convergence is strictly less than 1.  Moreover, they
found a relatively easily computable upper bound for the radius of
convergence in the case of a hyperbolic, formal space. They also
showed that if $A$ is a noetherian, connected graded commutative
algebra over a field $\mathbb F$ of characteristic zero and $\rho\sb A$ denotes
the radius of convergence of 
$$P\sb A(z)=\sum \sb {n\geq o}\dim_{\mathbb F} \operatorname{Tor}\sb n\sp A(\mathbb F,
\mathbb F)\cdot z\sp  n,$$
then either $\rho\sb A=+\infty$ and $A$ is a polynomial algebra; or
$\rho\sb A=1$, $A$ is a complete intersection, and the coefficients of
$P\sb A(z)$ grow polynomially; or $\rho\sb A<1$, $A$ is not a complete
intersection, and the coefficients of $P\sb A(z)$ grow
exponentially.  Avramov later generalized this result to any characteristic \cite {avramov2}.

The
written version of Avramov's Luminy talk on the close links between local algebra and rational homotopy
theory provides an excellent and thorough introduction to the subject \cite
{avramov3}.  His article contains the first ``dictionary" between
rational homotopy theory and local algebra, explaining how to translate
notions and techniques from one field to the other.  Given a
theorem in one field, applying the dictionary leads to a statement in the
other field that stands a reasonable chance of being true, though the
method of proof may be completely different.  

Avramov and Halperin wrote another thorough introduction to the subject  in the proceedings
of the Stockholm conference of 1983 \cite {avrhal}.  It begins at a more
elementary level than the survey article of Avramov in the proceedings of
the Luminy conference, leading the reader from first principles of
differential graded homological algebra to  notions of homotopy fiber and
loop space and on to the homotopy Lie algebra.  

In his introductory article \cite {avramov3}, Avramov emphasized the importance  of minimal models
in local ring theory.  If $K\sp R$ is the Koszul complex of a local ring, then there is a minimal, commutative
cochain algebra $(\Lambda V,d)$ over the residue field of $R$ that is 
quasi-isomorphic to $K\sp R$.  Avramov called $(\Lambda V,d)$ the
minimal model of $R$.  He established its relevance by observing that in
degrees greater than $1$, the graded Lie algebra derived from $(\Lambda
V,d\sb2)$ is isomorphic to the homotopy Lie algebra of $R$.

In \cite {fŽhalthom} F\'elix, Halperin, and Thomas continued the in-depth study of the
homotopy Lie algebra of a rationally hyperbolic space begun by F\'elix and
Halperin in \cite {fŽhal}.  They showed, for example, that if $E$ is rationally
hyperbolic, then its rational homotopy
Lie algebra is not solvable.   Moreover they proposed as
conjectures translations of their theorems into local algebra, where, for a
local ring $(R,\mathfrak m)$ with residue field $\mathbb F$, L.-S.
category is replaced by $\dim \sb {\mathbb F} (\mathfrak m/\mathfrak m\sp
2)-\operatorname {depth} R$ and infinite dimensional rational homotopy is
replaced by $R$ not being a complete intersection.  Recall that   
$$\operatorname {depth} R=\inf \{ j\mid \operatorname{Ext}\sp j\sb
R(\mathbb F,R)\not=0\}.$$

In \cite {fŽhalthom} F\'elix, Halperin and Thomas also mentioned a very important
conjecture due to Avramov and F\'elix, stating that the homotopy Lie
algebra of a rationally hyperbolic space should contain a free Lie algebra on
at least two generators.  This conjecture has motivated much interesting
work in the study of the homotopy Lie algebra and has not as yet (2006) been
proved.

Using minimal model techniques, Halperin and  B\o gvad proved two conjectures due to Roos, which are
``translations" of each other \cite{boghal}.  They  showed that
\begin{enumerate}
\item if $R$ is a local ring such that the Yoneda algebra
$\operatorname {Ext}\sb R\sp *(\mathbb F, \mathbb F)$ is noetherian, then $R$ is a
complete intersection; and
\item if $E$ is a simply connected, finite CW-complex such that the
Pontryagin algebra $H\sb *(\Omega E;\mathbb Q)$ is noetherian, then $E$ is
rationally elliptic.
\end{enumerate}
Their proof is based on a slightly weakened form of the Mapping Theorem
that holds over a field of any characteristic, as well as on ideas from the
article of F\'elix, Halperin and Thomas of the previous year \cite {fŽhalthom}.

In the spring of 1985 Halperin applied minimal model techniques to
answering on old question concerning the \emph{deviations} of a local
ring \cite {hal}.  The $j$th deviation, $e\sb j(R)$, of a noetherian, local, commutative ring
$R$ with residue field $\mathbb F$ is $\dim \sb {\mathbb F}\pi\sp j(R)$. 
Assmus had shown in 1959 that $R$ is a complete intersection if
and only if $e\sb j(R)=0$ for all $j>2$ \cite{assmus}, raising the question of
whether any deviation could vanish if $R$ were not a complete
intersection.   Halperin succeeded in answering this question, showing that
if $R$ is not a complete intersection, then
$e\sb j\not=0$ for all $j$.

The \emph{Five Author paper} \cite {5-auth} represents
a great leap forward in understanding the structure of the homotopy
Lie algebra of a  space or of a local ring.   The principal innovation of the
Five Author paper consists in exploiting the \emph{radical}
of the homotopy Lie algebra, i.e., the sum of all of its solvable ideals, which
rational homotopy theorists had begun to study in 1983.  The radical itself is
in general not solvable.  Recall, as mentioned above, that the rational homotopy Lie algebra of a rationally hyperbolic space is not solvable.

Expressed in topological terms, the main theorem of the Five Author
paper states that if $E$ is a simply connected CW-complex of finite type
and $\operatorname{cat}(E)=m<\infty$, then the radical of the
homotopy Lie algebra, $Rad(E)$ is finite dimensional and $\dim Rad(E)\sb
{\text{even}}\leq m$.  This is a consequence of two further theorems,
both of which are of great interest themselves.  The first concerns the
relations among the rational L.-S. category of a space and the depth and
global dimension of its homotopy Lie algebra.  Recall that the gobal
dimension of a local ring $R$ with residue field $\mathbb F$ satisfies
$$\operatorname{gl.dim.} (R)=\sup\{ j\mid \operatorname{Ext}\sb R\sp j(\mathbb F,
\mathbb F)\not=0\}.$$
The precise statement of this theorem in topological terms is then
that if $L$ is the homotopy Lie algebra of a simply-connected CW-complex
of finite type $E$, then either
$$\operatorname {depth}UL<\operatorname {cat}\sb 0(E)<\operatorname
{gl.dim.} UL$$
or
$$\operatorname {depth}UL=\operatorname {cat}\sb 0(E)=\operatorname
{gl.dim.} UL.$$
The second theorem states that under the same hypotheses, if
$\operatorname {depth}UL<\infty$, then $Rad(E)$ is finite dimensional and
satisfies $\dim Rad(E)\sb {\text{even}}\leq \operatorname {depth}UL$. 
Moreover, if $\dim Rad(E)\sb {\text{even}}=\operatorname {depth}UL$, then
$Rad(E)=L$.

Rational homotopy theorists have exploited extensively the results  of \cite {5-auth} in  developing a deep understanding of the homotopy Lie algebra of rationally
hyperbolic spaces.  The methods the five authors devised to prove their
results have turned out to be extremely important as well.   For example,
since their goal was to relate $\operatorname {cat}\sb 0(E)$ to
$\operatorname {depth}(L)$, they needed to construct a model of the
quotient cochain algebra $(\Lambda V/\Lambda \sp {>n} V,\bar d)$, where
$(\Lambda V,d)$ is the Sullivan minimal model of $E$.  Their method for
doing so, based on perturbation of a model for $(\Lambda V/\Lambda
\sp {>n} V,\overline {d\sb 2})$, proved to be useful in a number of other circumstances, such as in the proof of Theorem \ref{thm:mine}.



\begin{thebibliography}{99}

\bibitem {andrŽ} M. Andr\'e, \textsl {Cohomologie des alg\`ebres
diff\'erentielles o\`u op\`ere une alg\`ebre de Lie,} T\^ohoku Math. J. \textbf{14} (1962),  263--311.

\bibitem {anick1} D. J. Anick, \textsl {A counter-example to a conjecture of
Serre,} Annals of Math. \textbf{115} (1982),  1--33.

\bibitem {anick2} D. J. Anick, \textsl {Comment: ``A counter-example to a conjecture of
Serre,''} Annals of Math. \textbf{116} (1982),  661.

\bibitem  {assmus} E. Assmus, \textsl {On the homology of local rings,} Ill. J.
Math .\textbf{3} (1959),  187--199

\bibitem {avramov1} L. Avramov, \textsl {Homology of local flat extensions and
complete intersection defects,} Math. Ann.\textbf{ 228} (1977), 
27--37.

\bibitem{avramov2} L. Avramov, \textsl{Local rings of finite simplicial dimension,}  Bull. Amer. Math. Soc. \textbf{10} (1984) 289--291.

\bibitem {avramov3} L. Avramov, \textsl {Local algebra and rational
homotopy,} Homotopie Alg\'ebrique et Alg\`ebre Locale, 
Ast\'erisque, vol. 113--114, 1984, pp. 16--43. 



\bibitem {avrhal} L. Avramov and S. Halperin, \textsl {Through the looking
glass: a dictionary between rational homotopy theory and local
algebra,}  Algebra, Algebraic Topology and their
Interactions, Springer Lecture Notes in Mathematics, vol. 1183,
1986,  pp. 1--27. 



\bibitem {boghal} R. B\o gvad and S. Halperin, \textsl {On a conjecture of Roos,}
Algebra, Algebraic Topology and their Interactions,
Springer Lecture Notes in Mathematics, vol. 1183, 1986, 
120--126

\bibitem {bousgug} A. K. Bousfield and V. K. A. M. Gugenheim, \textsl {On PL De
Rham Theory and Rational Homotopy Type,} Memoirs of the A.~M.~S., vol. 179, 1976.

\bibitem {cartan}  E. Cartan, \textsl {Sur les nombres de Betti des espaces de groupes close,} C.R. Acad. Sci. Paris \textbf{ 187} (1928), 196--197.

\bibitem {dgms} P. Deligne, P. Griffiths, J. Morgan and D. Sullivan, \textsl {Real homotopy theory of K\"ahler manifolds} Inventiones \textbf{ 29} (1975),  245--274.

\bibitem {dwyspal} W. Dwyer and J. Spalinski, \textsl {Homotopy theories and model
categories}, Handbook of Algebraic Topology (I.~M.~James, ed.), 
North-Holland, 1995, pp.  73--126.

\bibitem {fŽhal} Y. F\'elix and S. Halperin, \textsl {Rational L.-S. category and
its applications,} Trans.  A.~M.~S. \textbf{ 273} (1982), 
1--37

\bibitem  {5-auth} Y. F\'elix, S. Halperin, C. Jacobsson, C. L\"ofwall, and J.-C.
Thomas, \textsl {The radical of the homotopy Lie algebra,} Amer. J.
Math. \textbf{110} (1988),  301--322.

\bibitem {fŽhallem} Y. F\'elix, S. Halperin and J.-M. Lemaire, \textsl {Rational
category and conelength of Poincar\'e complexes,} Topology \textbf{37} (1998),  743-748.


\bibitem {fŽhalthom} Y. F\'elix, S. Halperin and J.-C. Thomas, \textsl {The
homotopy Lie algebra for finite complexes,} Publ. Math. I.~H.~E.~S. \textbf{56} (1982),  387--410.

\bibitem {fht} Y. F\'elix, S. Halperin and J.-C. Thomas, \textsl{Rational
Homotopy Theory,} Graduate Texts in Mathematics, vol. 205, Springer-Verlag, 2001.

\bibitem {fŽthom} Y. F\'elix and J.-C. Thomas, \textsl {The radius of convergence
of Poincar\'e series of loop spaces,} Invent. Math. \textbf{ 68} (1982), 257--274.

\bibitem {goerss} P. Goerss, \textsl {Model categories and simplicial methods,} Lecture notes from the summer school on ``Interactions between Homotopy Theory and Algebra,'' Contemporary Math., 2006.

\bibitem {gugmunk}   V.K.A.M Gugenheim and H. J. Munkholm, \textsl {On the extended 
functoriality of Tor and Cotor,} J. Pure Appl. Algebra \textbf{ 4} (1974) ,  9--29.

\bibitem  {hal} S. Halperin, \textsl {The nonvanishing of deviations of a local ring,} Comment. Math. Helv.\textbf{ 62} (1987),  646--653.

\bibitem {hallem} S. Halperin and J.-M. Lemaire, \textsl {Notions of category in
differential algebra,} Algebraic Topology: Rational
Homotopy, Springer Lecture Notes in Mathematics, vol. 1318, 
pp. 138--153.


\bibitem {hess} K. Hess, \textsl {A history of rational homotopy theory,} History of Topology, Elsevier Science B.V., 1999,  pp. 757--796.



\bibitem {hovey} M. Hovey, \textsl{Model Categories}, Mathematical Surveys and Monographs, vol. 63, American
Mathematical Society, 1999.


\bibitem {iwase} N. Iwase, \textsl {Ganea's conjecture on
Lusternik-Schnirelmann category,} Bull.  London Math.  Soc. \textbf{30} (1998),  623--634.

\bibitem {levin} G. Levin, \textsl{Homology of local rings,} Ph. D.
Thesis, Univ. Chicago, 1965.

\bibitem{loday} J.-L. Loday, \textsl{Cyclic Homology (Second Edition)}, Grundlehren der Mathematischen Wissenschaften, vol. 30, Springer-Verlag, 1998.

\bibitem {loroos} C. L\"ofwall and J.-E. Roos, \textsl {Cohomologie des
alg\`ebres de Lie gradu\'ees et s\'eries de Poincar\'e-Betti non
rationnelles,} C. R. Acad. Sci. Paris \textbf{290} (1980) , 
733--736

\bibitem {mimo} J. Milnor and J. Moore, \textsl {On the structure of Hopf
algebras,} Annals of Math. \textbf{81} (1965),  211--264.

\bibitem {roos} J.-E. Roos, \textsl {Relations between the Poincar\'e-Betti
series of loop spaces and local rings,} S\'eminaire d'alg\`ebre
Paul Dubreil, Springer Lecture Notes in Mathematics, vol.
740, 1979, pp. 285--322

\bibitem {sjo} G. Sj\"odin, \textsl {Hopf algebras and derivations,} J.
Algebra \textbf{64} (1980),  218--229.

\bibitem {suvp} D. Sullivan and M. Vigu\'e-Poirrier, \textsl {The homology
theory of the closed geodesic problem,} J. Differential Geometry  \textbf{11} (1976),  633--644.

\end{thebibliography}
\end{document}